\documentclass[12pt]{amsart}

\usepackage{color}
\usepackage[centertags]{amsmath}
\usepackage[margin=1in]{geometry}
\usepackage{amsfonts}
\usepackage{amsthm}
\usepackage{newlfont}
\usepackage{verbatim}
\usepackage{amscd}
\usepackage{amsgen}
\usepackage{amssymb}
\usepackage{stmaryrd}
\usepackage{amssymb,amsmath}
\usepackage{amscd}
\usepackage{enumerate}
\usepackage{color}
\usepackage{xypic}
\usepackage{graphicx}

\newlength{\defbaselineskip} 
\newcounter{wsk}
\newcounter{wskk}
\theoremstyle{plain}
\newtheorem{thm}{Theorem}[section]
\newtheorem{cor}[thm]{Corollary}
\newtheorem{con}[thm]{Conjecture}
\newtheorem{df}[thm]{Definition}
\newtheorem{lema}[thm]{Lemma}
\newtheorem{obs}[thm]{Proposition}
\newtheorem{exm}[thm]{Example}
\newtheorem{question}[thm]{Question}

\newtheorem{tthm}[wsk]{Main Theorem}
\newtheorem{ghip}[wskk]{Conjecture}
\newtheorem{gcor}[wsk]{Corollary}

\newtheorem{rem}[thm]{Remark}
\newtheorem{pr}{Algorithm}
\theoremstyle{definition} 
\theoremstyle{definition} \newtheorem{ex}{Example}[section] %

 \numberwithin{equation}{section}
\def\p{\mathbb{P}}
\def\r{\mathbb{R}}
\def\z{\mathbb{Z}}
\def\n{\mathbb{N}}
\def\Z{\mathbb{Z}}
\def\N{\mathbb{N}}
\def\m{\mathfrak{m}}
\def\c{\mathbb{C}}

\def\nt{\mathfrak{nt}}
 \DeclareMathOperator{\Proj}{Proj}

\DeclareMathOperator{\rad}{Rad}

\DeclareMathOperator{\Spec}{Spec}
\def\p{\mathbb{P}}

\def\ob{\begin{obs}}
\def\kob{\end{obs}}
\def\dow{\begin{proof}}
\def\kdow{\end{proof}}
\def\kwadrat{\hfill$\square$}
\def\tw{\begin{thm}}
\def\ktw{\end{thm}}
\def\hip{\begin{con}}
\def\khip{\end{con}}
\def\lem{\begin{lema}}
\def\klem{\end{lema}}
\def\ex{\begin{exm}}
\def\prog{\begin{pr}}
\def\kprog{\end{pr}}
\def\wn{\begin{cor}}
\def\kwn{\end{cor}}
\def\uwa{\begin{rem}}
\def\kuwa{\end{rem}}
\def\kex{\end{exm}}
\def\dfi{\begin{df}}
\def\kdfi{\end{df}}
\def\flo{{\mathcal G}}
\def\soc{{\mathcal S}}
\definecolor{zielony}{rgb}{0.5, 0.9, 0.1}
\definecolor{czerwony}{rgb}{0.9, 0.2, 0.1}
\definecolor{niebieski}{rgb}{0.3, 0.1, 0.9}

\xyoption{line}
\begin{document}
\title{{Constructive degree bounds for group-based models}}
\author{Mateusz Micha\l ek}
\thanks{The scientific work is supported by Polish MNiSW grant 
Nr  IP2011 004971 in the years 2012-2014}
\maketitle
\begin{abstract}
Group-based models arise in algebraic statistics while studying evolution processes. They are represented by embedded toric algebraic varieties. Both from the theoretical and applied point of view one is interested in determining the ideals defining the varieties. Conjectural bounds on the degree in which these ideals are generated were given by Sturmfels and Sullivant \cite[Conjectures 29, 30]{SS}.
We prove that for the 3-Kimura model, corresponding to the group $G=\z_2\times\z_2$, the projective scheme can be defined by an ideal generated in degree $4$. In particular, it is enough to consider degree $4$ phylogenetic invariants to test if a given point belongs to the variety. We also investigate $G$-models, a generalization of abelian group-based models. For any $G$-model, we prove that there exists a constant $d$, such that for any tree, the associated projective scheme can be defined by an ideal generated in degree at most $d$.
\end{abstract}
\section{Introduction}
We investigate properties of a special class of not necessarily normal, projective toric varieties. Their construction is motivated by the applications of mathematics to phylogenetics. Such varieties are represented by a phylogenetic tree and a model of evolution. Each edge of a phylogenetic tree corresponds to a mutation. The probabilities of different mutations form a matrix, called the transition matrix. Biologists distinguished
certain types of matrices specified by a model of evolution. A given choice of transition matrices gives
us a probability distribution on the set of states of observed species. We may fix a model but vary entries
of transition matrices obtaining different probability distributions. This is an algebraic map. The
closure of its image is a variety that is the main object of our study.

We will not describe the connections with applications, referring to \cite{Ph, PachSturm}. On the other hand we give an exact mathematical definition of the class of varieties we consider. 
Let $T$ be a rooted tree, that is a connected graph without cycles and with one distinguished vertex $v_0$ called the root. We direct the edges of $T$ away from the root. Let $V$ be the set of vertices of $T$. For each $v\in V$, consider a vector space $W_v$ with a fixed basis $e_1^v,\dots,e_{n(v)}^v$. We will construct varieties embedded in $\bigotimes_{v\in L} W_v$, where $L$ is the set of leaves of $T$. The variety is given as the closure of a parametrization map
$$\psi:W_{v_0}\times\prod_{(v_1,v_2)}W_{v_1}^*\otimes W_{v_2}\rightarrow \bigotimes_{v\in L} W_v,$$
where the product $\prod_{(v_1,v_2)}$ is taken over all edges of $T$. Consider the basis $\bigotimes_{v\in L} e_{j(v)}^v$ of the vector space $\bigotimes_{v\in L} W_v,$ indexed by functions $j$ that associate to a leaf $v$ an index from $1$ to $n(v)$. We define $\psi$ by:
$$(\bigotimes_{v\in L} e_{j(v)}^v)^*(\psi((w_0,\prod_{(v_1,v_2)} w_{v_1^*,v_2})))=\sum_i (e_{i(v_0)}^{v_0})^*(w_0)\prod_{(v_1,v_2)}(e_{i(v_1)}^{v_1}\otimes e_{i(v_2)}^{v_2*})(w_{v_1^*,v_2}),$$
where $w_{v_1^*,v_2}\in W_{v_1}^*\otimes W_{v_2}$ and the sum $\sum_i$ is taken over all functions that associate to a vertex $v$ an index of a basis vector of $W_v$ and agree with $j$, when $v$ is a leaf. As usually by $^*$ we denote the dual vector, with respect to the chosen basis. Note that the definition above is just a simple probability computation according to a Markov process on a tree. Indeed, the function $j$ fixes basis vectors at leaves - this corresponds to fixing states. The parameter $w_0\in W_{v_0}$ gives the root distribution. Each vector in $W_{v_1}^*\otimes W_{v_2}$ can be identified with a transition matrix. The summation over functions $i$ corresponds to the summation over all states of vertices that extend given states of leaves. The formula follows from the Markov property.

The image of $\psi$ is the variety associated to \emph{the general Markov model}. Let $K_{1,n}$ be the \emph{claw tree} -- cf. Definition \ref{clawtree}, with the inner vertex $v_0$. We encourage the reader to check that for this tree and $W_{v_0}$ one dimensional, the obtained variety is the Segre embedding of spaces $W_1,\dots, W_n$ associated to leaves. Moreover, when $W_{v_0}$ is $k$ dimensional, the Zariski closure of the image of $\psi$ is the $(k{-}1)$-st secant variety of the Segre variety. More details on this construction, including motivation, can be found in \cite[Chapter 4]{PachSturm} or a short paper \cite{4aut}.

We are interested in equivariant versions of the above construction. First examples of such appeared in phylogenetics, the science that aims at reconstructing the evolution of DNA \cite{2Kimura, 3Kimura, Jukes.Cantor1969}. It was observed that certain symmetries between nucleobasis of DNA can be induced by an action of a group. Important mathematical implications of such a description of the variety were proven in \cite{HendyPenny}\footnote{We would like to thank Elizabeth Allman for bringing this fact to our attention.}.  

The models that are the central object of study of this article are called \emph{group-based}. 
In this setting we assume that all vector spaces $W_{v}\cong W$ are regular representations of a fixed finite abelian group $G$. Thus we can identify the basis with group elements. We restrict the domain of $\psi$ to $G$-invariant vectors. More precisely we consider:
$$\psi_G: W^G\times\prod_{(v_1,v_2)}(W^*\otimes W)^G\rightarrow \bigotimes_{v\in L} W,$$
where the action of the group on the tensor product and on the dual is induced from the action on the vector space.
Note that $W^G$ is one dimensional, hence the image of $\psi_G$ does not change after restricting the domain to $\prod_{(v_1,v_2)}(W^*\otimes W)^G$. The fact that $W^G$ is one dimensional corresponds in biology to the assumption of \emph{uniform root distribution}. The Zariski closure of the image of $\psi_G$ will be denoted by $X(T,G)$ and is known to be toric -- for a modern reference we advise the reader to consult \cite{SS}.
 We also let $X_\p(T,G)$ be the projective variety, such that $X(T,G)$ is the affine cone over $X_\p(T,G)$.

More generally, one can consider the vector spaces $W_v$ to be arbitrary representations. The model in such a case is called \emph{equivariant} -- it was introduced in \cite{DK}. This construction does not always lead to a toric variety. Indeed, the first presented example -- secant varieties of Segre varieties -- is an equivariant model with the trivial group.

The study of equivariant models was reduced to the case of claw trees. Still, the description of the varieties associated to claw trees is difficult. Let us present what is known.
\begin{enumerate}
\item The equivariant model associated to the trivial group corresponds to secant varieties of the Segre embedding. These varieties are an object of intensive studies. The description of the first secant was recently given in \cite{Raicu}, solving the GSS conjecture. It is known that an analogous description, by flattenings, for higher secants does not hold. On the other hand, it was proven that, for a fixed $k$, there exists $d$, such that the $k$-th secant of the Segre variety, for any number of factors, is set theoretically described by equations of degree at most $d$, \cite{DK2}.
\item The group-based model with $G=\z_2$ is the only model where we can explicitly give a description of the ideal \cite[Theorem 28]{SS}, \cite{Sonja}. It is definitely the simplest model one can consider. However, it turned out to have very interesting properties. One of them in the theorem of Buczy\'nska and Wi\'sniewski, which states that varieties associated to trivalent trees with the same number of leaves belong to one flat family \cite[Theorem 2.24]{BW}. Moreover the model has connections with different branches of mathematics \cite{SX, Man}. For the generalizations of the construction to graphs we advise the reader to consult \cite{BBKM}.
\item Equivariant models for an abelian group $G$ include all previous examples. Recently, Draisma and Eggermont proved that set theoretically they can be defined in (some) bounded degree \cite{DE}.
\end{enumerate}

There are several reasons why mathematicians are interested in equations defining the variety associated to a model. From the point of view of applications in phylogenetics one is interested in determining if a given point $Q$, representing a probability distribution, belongs to the image of $\psi$. One of the possible methods is to evaluate polynomials vanishing on the image on $Q$. These polynomials -- elements of the ideal defining $X(T,G)$, are called phylogenetic invariants. They are sought by people dealing with phylogenetics. Due to numerical evidence for small trees Sturmfels and Sullivant posed the following conjecture.
\begin{ghip}[Conjecture 30 \cite{SS}]\label{Kimura}
The ideal for the 3-Kimura model, that is a group-based model for $G=\z_2\times\z_2$, is generated in degree at most $4$ for any tree $T$.
\end{ghip}
One is particulary interested in this model, as it was introduced by theoretical biologists \cite{3Kimura}. Four elements of the group correspond to four nucleobasis forming the DNA and the group action captures symmetries between them.
\begin{tthm}\label{glKimura}
For any tree $T$ the \emph{projective} scheme $X_\p(T,\z_2\times\z_2)$ can be defined by an ideal generated in degree at most $4$.
\end{tthm}
Conjecture \ref{Kimura} has the following generalization, also due to Sturmfels and Sullivant.
\begin{ghip}[Conjecture 29 \cite{SS}]\label{ogolna}
For any group $G$ the ideal of the group-based model is generated in degree at most $|G|$ for any tree $T$.
\end{ghip}
A special class of equivariant models that generalizes abelian group-based models are $G$-models. In this setting we assume that $G$ is a group (not necessary abelian) and $W$ is a representation of $G$. Moreover, we assume that $G$ contains a normal, abelian subgroup $H$ such that $W$ is a regular representation of $H$. Thus a $G$-model is in fact a submodel of a group-based model corresponding to $H$. These models were introduced and studied in \cite{mateusz}.
An example of a $G$-model model that is a strict submodel of a group-based model is the 2-Kimura model and the Neyman model, also known as Jukes-Cantor model, for at least $3$ states.

\begin{tthm}\label{glgen}
For any $G$-model there exists a constant $d$ such that the associated projective scheme $X_\p$ can be defined by an ideal generated in degree at most $d$ for any tree $T$.
\end{tthm}
 The result presented above is quite similar to the main results of \cite{DE}, that generalize those in \cite{DK}. Let us state $3$ major differences:
\begin{enumerate}
\item The results in \cite{DE} concern a wide class of \emph{abelian equivariant} models, while ours deal with \emph{$G$-models} -- in particular, the general Markov model (resp. the 2-Kimura model) belongs to the first (resp. second) class and does not belong to the second (resp. first);
\item Our results are \emph{scheme theoretic}, while the results in \cite{DE} are \emph{set theoretic};
\item Our results are \emph{constructive}, while methods of \cite{DE} are \emph{existential}.
\end{enumerate}
In \cite{DE} the authors use a beautiful technique that is based on noetherian arguments for rings that are not finitely generated, but are equipped with an additional monoid action. For more on this technique see \cite{Sullfinite, Abraham, Dfinit}. Our methods are typical for toric geometry. Indeed, we believe that the presented approach can be used for much broader class of problems in toric geometry to determine the degree of generation -- see for example \cite{Bruns}. In particular, the methods cover group-based models for small groups $G$. As our proofs are constructive, for group-based models, one can explicitly write a polynomial algorithm for testing if a point belongs to $X(T,G)$. The existence of such an algorithm was proven in \cite{DE}. Still, our bounds on the degree in which the schemes are defined are weak. This implies that the complexity of the algorithm would be in general too big to be of practical value. However, for the 3-Kimura model our results on the degree are sharp.

\begin{gcor}\label{wniosekwstep}
For any tree $T$, a point $Q$ belongs to $X(T,\z_2\times\z_2)$ if and only if all phylogenetic invariants of degree at most four vanish at $Q$.
\end{gcor}
We believe that Corollary \ref{wniosekwstep} may be used in applied mathematics. Indeed, we know that for 3-Kimura model it is enough to consider phylogenetic invariants of degree at most $4$. To generate all invariants of a given degree one can apply for example methods of \cite[Section 3.2]{DBM}. Hence, one can obtain enough phylogenetic invariants to test if a point belongs to a variety associated to any tree.

In Section \ref{construction} we present a description of polytopes associated to varieties $X(T,G)$.
The general method that we use can be applied to other problems in toric geometry. It is described in Section \ref{idea}. In Section \ref{combinatoriallemma} we present a basic combinatorial lemma. It is the main tool to prove Main Theorem \ref{glgen}, which we do in Section \ref{proof2}. The reader may notice the resemblance between \emph{contractions} of tensors and methods of Section \ref{proof2}.
The technical proof of Main Theorem \ref{glKimura} can be found in Section \ref{proof}. In the last section \ref{open} we present open problems in terms of algebraic geometry.

\section*{Acknowledgements}
The article was started at Universite Grenoble I and finished during a visit at Max Planck Institute in Bonn. I am grateful for their hospitality, especially to Laurent Manivel. I would like to thank Jaros\l aw Wi\'sniewski for introducing me to the subject. I appreciate important remarks of two anonymous reviewers.
\section{Notation}
As the paper presents two different main theorems, related to the same topic, some notation conventions are settled. The most important objects are: a finite group $G$, a tree $T$, a group of group-based flows $\flo$ (Definition \ref{flows}), polytope $P$. In each section additional objects are introduced. In general, objects defined by small latin letters are local and their meaning may differ from one section to another. Throughout the paper the group-based flows are in constant use. Group-based flows are specific functions with an additive group structure. They are often denoted by small latin letters (possibly with a subscript). However, as the reader may find out in Section \ref{construction}, one can naturally identify group-based flows with integral points of a polytope $P$. Thus, when representing them as integral points, it is more natural to denote them by capital latin letters. Moreover, the addition operation changes. Adding group-based flows as functions, gives a group-based flow. However, adding two integral points corresponding to group-based flows, gives another integral point, that may not (and indeed never will) represent a group-based flow.

By a lattice we always mean integral lattice, that is a group isomorphic to $\z^n$ for some $n\in\N$. We say that a homogeneous ideal $I$ is generated in degree $d$ if all its elements of degree less or equal to $d$ generate the ideal.
\section{Constructions}\label{construction}
Let $T$ be a rooted tree and $G$ a finite abelian group. We assume that all edges of $T$ are directed away from the root, however we could consider any orientation.
\dfi[Edges $E$, Vertices $V$, Leaves $L$, Node $N$]
We define the sets $E$ and $V$ to be respectively the set of edges and vertices of the tree $T$. Vertices of degree one are called leaves. The set of leaves is denoted by $L$. By abuse of notation, edges adjacent to leaves will also be called leaves. We define the set of nodes $N:=V\setminus L$.
\kdfi

There are three canonical groups one can consider.
\dfi[Edge labellings $G_E$, Vertex labellings $G_V$, Nodes labellings $G_N$]
The group $G_E$ of edge labellings consists of all functions $f:E\rightarrow G$, with the group operation defined by $(f_1+f_2)(e)=f_1(e)+f_2(e)$. The group $G_E$ is isomorphic to $G^{|E|}=G\times\dots\times G$.

The group $G_V$ of vertex labellings consists of all functions $g:V\rightarrow G$, with the group operation defined by $(g_1+g_2)(v)=g_1(v)+g_2(v)$. The group $G_V$ is isomorphic to $G^{|V|}$.

The group $G_N$ of node labellings consists of all functions $g:N\rightarrow G$, with the group operation defined by $(g_1+g_2)(v)=g_1(v)+g_2(v)$. The group $G_N$ is isomorphic to $G^{|N|}$.
\kdfi
As the tree $T$ is directed there is a canonical morphism $s:G_E \rightarrow G_V$.
\dfi[The summing morphism $s$]
We define $s:G_E\rightarrow G_V$ by
\begin{equation}\label{eqnsum}
s(f)(v)=\sum_{e=(x,v)}f(e)-\sum_{e=(v,x)}f(e),
\end{equation}
where the first sum is taken over all edges incoming to $v$ and the second over outgoing edges.
\kdfi

In phylogenetics it is natural to distinguish the set of leaves and nodes.
\dfi[Projection $\pi_N$]
We define the projection $\pi_N: G_V\rightarrow G_N$ by restricting the domain of a function from $V$ to $N$.
\kdfi
The most important combinatorial objects that we consider are the following.
\dfi[Group-based flows]\label{flows}
We define the group $\flo$ to be the kernel of the composition $\pi_N\circ s:G_E\rightarrow G_N$. We call elements of $\flo$ group-based flows. We have got an exact sequence:
$$0\rightarrow\flo\rightarrow G_E\rightarrow G_N\rightarrow 0.$$
As an abstract group we have an isomorphism $\flo\cong G^{|L|-1}$. The group $\flo$, up to isomorphism, is independent on the choice of orientation of edges.
\kdfi
\uwa
\emph{
The elements of $\flo$ are such edge labellings that "the signed sum around each node is trivial". If we associated to edges elements of $(\r,+)$ instead of $G$, then we would obtain a well-known condition for a flow.
Notice that we would allow associations of negative numbers to edges. In particular we would not distinguish sources and sinks -- both would be just vertices $v$ at which the sum $s(f)(v)$, defined in \eqref{eqnsum}, can be nonzero. Thus sources and sinks would be leaves.}

\emph{
This analogy gives a justification for the name "group-based flows". The group of group-based flows is important not only in phylogenetics. As it was observed by Manon , group-based flows appear while studying conformal field theory \cite{Man, Manonneww}.
}

\emph{Note also that the construction of group-based flows can easily be generalized to other subsets of vertices by considering arbitrary $N'\subset V$ and defining $\pi_{N'}$ respectively.}
\kuwa

\ex
\emph{
Let us consider the group $G=\z_3$ and the following tree rooted at the top vertex:
$$
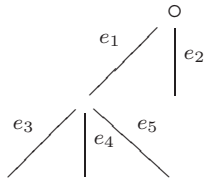
\xymatrix{
&&\circ\ar@{-}[dl]_{e_1}\ar@{-}[d]^{e_2}\\
&\ar@{-}[dl]_{e_3}\ar@{-}[d]^{e_4}\ar@{-}[dr]^{e_5}&\\
&&\\
}$$
Here $e_2$, $e_3$, $e_4$ and $e_5$ are leaves. An example of a group-based flow is an association $e_1\rightarrow 2$, $e_2\rightarrow 1$, $e_3\rightarrow 1$, $e_4\rightarrow 2$, $e_5\rightarrow 2$.}
\kex
\uwa\emph{
Eventually we will see in Definitions \ref{P} and \ref{X}, that group-based flows correspond to vertices of the polytope associated to the toric variety. Such objects were first introduced only for the group $\z_2$ and trivalent trees in \cite[Definition 3.1]{BW} and were called networks. They were originally defined as pairwise disjoint paths, beginning and ending at leaves. Definition \ref{flows} is a direct generalization. Let us define a bijection between the set of networks and the set of group-based flows for $G=\z_2$. Fix a network $N$. Let $S$ be the set of all edges comprised in any path belonging to the network $N$. We construct a group-based flow as follows:}
$$f(e)=
\begin{cases}
1\text{ for }e\in S,\\
0\text{ for }e\not\in S.
\end{cases}
$$
\emph{
For the group $\z_2$ both constructions do not depend on the orientation of the tree. Below we give an example of the bijection on a specific element. On the left we present a network, where the edges belonging to paths are not dashed. On the right there is the corresponding group-based flow.}
$$
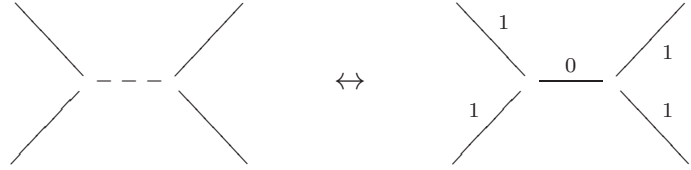
\xymatrix{
\ar@{-}[dr]&&&\ar@{-}[dl]&&\ar@{-}[dr]^1&&&\ar@{-}[dl]^1\\
&\ar@{--}[r]&\ar@{-}[dr]&&\leftrightarrow&&\ar@{-}[r]^0&\ar@{-}[dr]^1&\\
\ar@{-}[ur]&&&&&\ar@{-}[ur]^1&&&\\
}$$
\kuwa
It is well-known that, if the graph is a tree, it is enough to know the values at sources and sinks to reconstruct the flow. Of course, one cannot take arbitrary values for sources and sinks. For any flow, the sum over sources and sinks has to be equal to zero.
\dfi[Sockets $\soc$]\label{socket}
We define the group of sockets $\soc$ whose elements are functions $f:L\rightarrow G$ such that $\sum_{l\in L}f(l)=e_0$, where $e_0$ is the neutral element of the group. The group operation is defined by $(f_1+f_2)(e)=f_1(e)+f_2(e)$.
\kdfi
Notice that by restricting a group-based flow to leaves one obtains a socket.
In fact one can easily prove the following.
\ob
The groups $\flo$ and $\soc$ are naturally isomorphic.
\kwadrat
\kob
\uwa\emph{
As it was for group-based flows also sockets were first defined in \cite[Definition 3.1]{BW} for trivalent trees and the group $\z_2$. They were defined as even subsets of leaves. We see that the condition that the subsets are even  corresponds to the summation condition in Definition \ref{socket}.}
\kuwa
We come to the definition of the polytope $P$ that represents the variety $X_\p(T,G)$. 
\dfi[Lattices $M_e$, $M_E$, Basis $b_{(e,g)}$]
For a fixed edge $e$ we define $M_e$ as the lattice with basis elements indexed by all pairs $(e,g)$ for  $g\in G$. The basis element indexed by $(e,g)$ is denoted by $b_{(e,g)}\in M_e$. We define $M_E=\prod_{e\in E} M_e$. The elements $b_{(e,g)}$ form its basis for $(e,g)\in E\times G$. In particular $M_E\simeq \z^{|E|\cdot|G|}$.
\kdfi
\dfi[Polytope $P$]\label{P}
To a group-based flow $f\in\flo$ one can naturally associate an element $P_f:=\sum_{e\in E} b_{(e,f(e))}\in M_E$. This is an element whose coordinates are equal to either $1$ or $0$.

We define $P$ to be the convex hull of all the points $P_f$ over all group-based flows $f$.
\kdfi
The following proposition follows from the fact that $P$ is a subpolytope of the unit cube in $M_E$.
\ob
All integral points of $P$ are vertices and are of the form $P_f$ for certain group-based flow $f\in\flo$.\kwadrat
\kob
\dfi[Monoid $S(P)$, Semigroup algebra $\c{[}P{]}$]\label{monoidalgebra}
We define the monoid $S(P)$ as the submonoid of $M_E$ spanned by vertices of $P$. We define $\c[P]$ as the semigroup algebra associated to the monoid $S(P)$. Precisely $\c[P]$ is a complex vector space with basis elements identified with elements of $S(P)$. The multiplicative structure on $\c[P]$ is induced from the additive structure on $S(P)$.
\kdfi
\dfi[Variety $X(T,G)$]\label{X}
Notice that $\c[P]$ is a graded algebra, taking vertices of $P$ to be of degree $1$. We define the projective variety $X_\p(T,G):=\Proj \c[P]$ and $X(T,G):=\Spec \c[P]$.
\kdfi
It is not obvious that the variety $X(T,G)$ agrees with the one defined in the Introduction. Indeed, it is a nontrivial theorem -- \cite[p. 346/347]{SS, mateusz}.
\uwa\emph{
The variety $X_\p(T,G)$ is \emph{not} the toric variety associated to the polytope $P$ as defined for example in \cite[Section 1.5]{Fult}. The variety $X_\p(T,G)$ does not have to be normal -- for example when $G=\z_6$ \cite[Computation 4.3]{DBM}. In fact, the toric variety associated to the polytope $P$ is the normalization of $X_\p(T,G)$.
}
\kuwa

Consider the vector space $V'$ with basis elements indexed by vertices of $P$. There is a natural embedding of $X_\p(T,G)$ in the projective space $\p(V')$. The choice of the basis induces a toric structure on $\p(V')$. The projective space $\p(V')$ contains a dense algebraic torus $\mathbb{T}$. The intersection of $\mathbb{T}$ with $X_\p(T,G)$ is also an algebraic torus that acts, with a dense orbit, on $X_\p(T,G)$. We are interested in the ideal $I(T,G)$ that defines $X_\p(T,G)$ in $\p(V')$.

The following crucial fact is well known in toric geometry.
\ob[\cite{Stks} Lemma 4.1]
Let $A$ be any finite set in a lattice $M$. Define the corresponding ideal $I_A$ as in \cite[Chapter 4]{Stks}.
The toric ideal $I_A$ is spanned, as a vector space, by the set of binomials
$$\prod x_i^{a_i}-\prod x_i^{b_i},$$
where $\sum a_iA_i=\sum b_iB_i$ encodes an integral relation between points $A_i,B_i\in A$.
\kwadrat
\kob
We apply this proposition to $A=P$.
We see that the description of the ideal $I(T,G)$ is \emph{reduced to the study of integral relations between the vertices of the polytope P}.

Let us briefly introduce $G$-models. We assume that a finite group $G$ contains a normal, abelian subgroup $H$. In particular we can define the polytope $P$ corresponding to the group-based model for $H$. Consider the set $O$ of orbits of the adjoint action of $G$ on $H$.
\dfi[$M_{e,sub}$, $M_{E,sub},P_{sub}$]
Fix an edge $e\in E$. We define $M_{e,sub}$ as the lattice with basis elements indexed by all pairs $(e,o)$ for  $o\in O$. The basis element indexed by $(e,o)$ is denoted by $b_{(e,o)}\in M_{e,sub}$. We define $M_{E,sub}:=\prod_{e\in E} M_{e,sub}$. There is a natural morphism $M_E\rightarrow M_{E,sub}$, where $M_E$ is the lattice for the group-based model associated to $H$. The morphism associates to the basis vector $b_{(e,g)}$ the basis vector $b_{(e,[g])}$, where $[g]$ is the orbit of $g$. Let $P_{sub}$ be the image of $P$. The algebraic variety $X$ associated to the $G$-model is $\Spec \c[P_{sub}]$.
\kdfi

\subsection{Group-based polytopes}
So far we have seen that the group $\flo$ and polytope $P$ are strongly related. We believe that the object consisting of such a pair deserves a separate, general definition.
\dfi[Group-based polytope]\label{grbasedpoly}
Let $Q$ be a lattice polytope in the lattice $M$. Suppose that a group $G$ acts on the integral points of $Q$. We say that the pair $(Q,G)$ is a group-based polytope if the action of $G$ preserves the linear relations. That is, for any relation $\sum q_i=\sum q_j$ and any $g\in G$ we have $\sum g(q_i)=\sum g(q_j)$, where $q_i,q_j\in Q$ are lattice points.
\kdfi
As before, the polytope $Q$ generates a monoid $S(Q)$ which defines the semigroup algebra $\c[Q]$ -cf. Definition \ref{monoidalgebra}. The toric variety $\Spec \c[Q]$ has a natural embedding in the affine space $V'$, with basis elements corresponding to lattice points of $Q$. As $G$ acts on $Q$ we see that $V'$ is a representation of $G$. The assertion in Definition \ref{grbasedpoly} that the $G$ action preserves the linear relation is equivalent to the fact that the action restricts to the variety $\Spec \c[Q]$.

It is straightforward to see that $(P,\flo)$ is a group-based polytope. Moreover, as the action of $\flo$ on vertices of $P$ is transitive and free we obtain the following proposition.
\ob\label{action}
The natural ambient space $V'$ of the variety $X(T,G)$ is the regular representation of the group $\flo$. The induced action on $\p(V')$ restricts to the variety $X_\p(T,G)$.
\kwadrat
\kob
\uwa\emph{
In fact, canonically one should define group-based flows as labellings of edges by \emph{characters} of $H$, that is elements of $H^*$. As the groups $H$ and $H^*$ are isomorphic, in all constructions we get isomorphic objects. We decided not to introduce $H^*$ to avoid unnecessary confusion.}
\kuwa
\section{General methods in phylogenetics}
Although the construction of the variety $X_\p(T,G)$ takes several steps, still all the data that was used was the tree $T$ and the group $G$. Thus one would expect to find certain relations between invariants of the ideal $I(G,T)$, the group $G$ and the tree $T$.
\uwa\emph{
 An example of such a relation, in a slightly different setting, was presented in \cite{BBKM}. The genus of the graph plus one was proved to be the upper bound on the degree of generation of an algebra associated to any graph and the group $\z_2$.}
\kuwa
\dfi[Claw tree $K_{1,n}$]\label{clawtree}
The claw tree $K_{1,n}$ is the tree with one node and $n$ leaves.
\kdfi
We introduce certain invariants that are of great interest \cite[Chapter 5]{SS}.
\dfi[$\phi(G,T), \phi(G), \phi(G,n)$]
Let $\phi(G,T)$ be the maximal degree of the minimal set of generators of $I(G,T)$. Let $\phi(G,n):=\phi(G,K_{1,n})$. Let $\phi(G)=\sup \phi(G,T)$, where the supremum is taken over all trees.
\kdfi
Let us restate the conjecture of Sturmfels and Sullivant.
\hip[\cite{SS}, Conjecture 29]\label{generalconj}
We have $\phi(G)\leq |G|$.
\khip
It was separately stated in the most interesting case for $G=\z_2\times\z_2$, \cite[Conjecture 30]{SS}.
Note that, a priori, we do not know if $\phi(G)$ is finite, and indeed it is an open problem for groups different from $\z_2$.

Let us present the combinatorial way of representing group-based flows and relations between vertices of $P$. Using the isomorphism $G_E\cong G^{|E|}$, each element of $G_E$ can be represented as a sequence of group elements of length $|E|$. We represent group-based flows, as such sequences, presented as column vectors. A multiset of group-based flows is represented by a matrix. Each column corresponds to a group-based flow and the rows are indexed by edges of the tree $T$. Recall that group-based flows correspond to vertices of $P$. Consider any combination $\sum a_iv_i$ where $v_i\in P$ are vertices of $P$ and $a_i\in \n$. Such a combination can be encoded as a matrix, where each group-based flow corresponding to $v_i$ appears $a_i$ times. Of course such matrices should be considered up to permutation of columns.

Continuing this, we may encode a relation $\sum a_iv_i=\sum a_j v_j$ as a pair of matrices $A_1,A_2$ of the same size such that:
\begin{enumerate}
\item rows are indexed by edges,
\item each column represents a group-based flow,
\item for each edge $e$ the $e$-th row in $A_1$ is equal to the $e$-th row of $A_2$ up to permutation.
\end{enumerate}
\ex\emph{
Consider the following tree and numbering of edges:
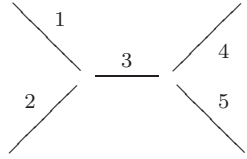
\begin{equation}\label{drzewo}\end{equation}
$$\xymatrix{
\ar@{-}[dr]^1&&&\ar@{-}[dl]^4\\
&\ar@{-}[r]^3&\ar@{-}[dr]^5&\\
\ar@{-}[ur]^2&&&\\
}$$
Take $G=\z_2$.
An example of a nontrivial relation between group-based flows is as follows:}
\[
\left[
\begin{array}{cccccccc}
1&0\\
1&0\\
0&0\\
0&1\\
0&1\\
\end{array}
\right],
\left[
\begin{array}{cccccccc}
1&0\\
1&0\\
0&0\\
1&0\\
1&0\\
\end{array}
\right].
\]
\kex
One of the main methods in algebraic phylogenetics is reduction to simpler trees. Each tree $T$, different from the claw tree, can be subdivided into two simpler trees $T_1$ and $T_2$, by dividing an inner edge.

$$
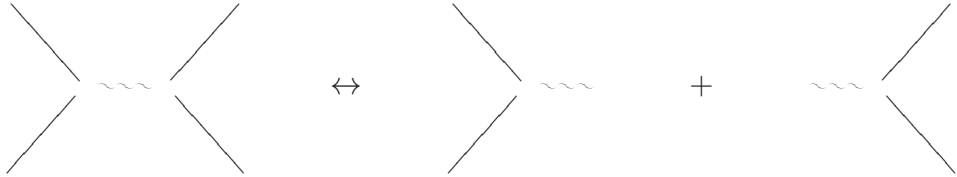
\xymatrix{
\ar@{-}[dr]&&&\ar@{-}[dl]&&\ar@{-}[dr]&&&&&&\ar@{-}[dl]\\
&\ar@{~}[r]&\ar@{-}[dr]&&\leftrightarrow&&\ar@{~}[r]&&+&\ar@{~}[r]&\ar@{-}[dr]&\\
\ar@{-}[ur]&&&&&\ar@{-}[ur]&&&&&&\\
}$$

The operation can be reversed. If we have got two trees with distinguished leaves we can join them, obtaining one tree.

Often one can transfer certain results from the trees $T_1$ and $T_2$ to the tree $T$. This is the case for the degree of generation.
We encourage the reader to prove the following result using the language of group-based flows.
\ob
Suppose that we have got two trees $T_1$ and $T_2$ with distinguished leaves and that a tree $T$ is obtained by joining them. Then $\phi(G,T)\leq \max(\phi(G,T_1),\phi(G,T_2))$.\kwadrat
\kob
In fact one can derive the generators of the ideal $I(G,T)$ from the ideal $I(G,T_1)$ and $I(G,T_2)$, \cite[Theorem 26]{SS}, \cite[Corollary 2.11]{Sethtfp}. Thus if we want to bound the degree in which $I(G,T)$ is generated for any tree $T$ it is enough to consider claw trees.
\uwa\emph{
The "reduction to claw trees" has been obtained in a much more general setting than the presented in this paper \cite{AllRhMarkov, DK}. Still the problem of describing the ideal for claw trees is extremely hard \cite[p.~637]{DK}, apart from the simple case $G=\z_2$ \cite{Sonja}.}
\kuwa

In algebraic geometry varieties can be compared on different levels. The most common distinction is to compare them as \emph{schemes} or as \emph{sets}. Consider two ideals $I_1,I_2$ in a ring of polynomials $\c[x_1,\dots,x_n]$. It is well-known that they define the same affine scheme if and only if they are equal. By Hilbert's Nullstellensatz we know that they define the same set if and only if their radicals are equal. We call the first equality \emph{ideal theoretic} and the second one \emph{set theoretic}. If the ideals are homogeneous they define not only affine schemes, but also projective schemes. Recall that $I:J^{\infty}=\{f: fJ^N\subset I \text{ for some }N\text{ large enough}\}$ is called the saturation of $I$ with respect to $J$. Let $\m=(x_1,\dots,x_n)$ be the irrelevant ideal. The homogeneous ideals $I_1$ and $I_2$ define the same projective scheme if and only if their saturations with respect to the irrelevant ideal are equal, that is $I_1:\m^\infty=I_2:\m^\infty$. We call equality of ideals after saturation \emph{scheme theoretic}.
Of course equality of projective schemes implies equality of projective sets. So if both $I_1$ and $I_2$ are simultaneously either contained or not in $\m$ then \emph{scheme theoretic} equality implies \emph{set theoretic}.
\ex\emph{
Consider $\c[x,y]$ and the ideals $I_1=(x,y)$, $I_2=(x^2,xy,y^2)$, $I_3=(1)$. All of them define the same (empty) projective scheme. Their saturation is $I_3$. The first two define a nonempty affine set, supported at $0$, with different scheme structure.}
\kex

\section{Idea}\label{idea}
\dfi[Subideal $I_d$]
For any graded ideal $I$ we define the subideal $I_d$ as the ideal generated by all elements in $I$ of degree at most $d$.
\kdfi
Let us start from general remarks on toric geometry. Consider a lattice polytope $Q$ contained in a hyperplane not passing through $0$. This allows us to put a grading to the associated monoid generated by $Q$. Hence the associated toric ideal $I_Q$ is homogeneous. Suppose that we want to prove that $I_Q$ and $I_{Qd}$ for some $d\in \n$ define the same projective scheme. Let $X$ be the projective variety defined by $Q$.
The integral points of $Q$ correspond to coordinates of the ambient projective space of $X$. Proving that the saturation of $I_{Qd}$ equals $I_Q$ is equivalent to proving that $I_{Qd}$ and $I_Q$ are equal in each localization with respect to any coordinate, represented by a lattice point $R\in Q$. Thus we have to prove that any generator of $I_Q$ multiplied by a sufficiently high power of the variable corresponding to $R$ belongs to $I_{Qd}$.

Let us restate this condition in combinatorial terms. The generators of $I_Q$ correspond to relations between points of $Q$. Let us fix a relation $\sum A_i=\sum B_j$, where $A_i,B_j\in Q$. Multiplying the corresponding element of the ideal by the variable corresponding to $R$ is equivalent to adding $R$ to both sides of the relation. Thus we have to prove that the binomial corresponding to the relation $\sum A_i+mR=\sum B_j+mR$ is generated by binomials from $I_Q$ of degree at most $d$ for $m$ sufficiently large.

A binomial corresponding to a relation $\sum R_i=\sum S_i$ between points of a polytope is generated in degree $d$ if and only if one can transform $\sum R_i$ to $\sum S_i$ using a sequence of following transformations. In each single transformation one can replace points $R_1,\dots,R_k$ for $k\leq d$ by $R_1',\dots,R_k'$ if they satisfy the relation $\sum_{i=1}^k R_i=\sum_{i=1}^k R_i'$. In such a case we say that the relation is generated in degree $d$.

The strategy of the proof is very simple:
\begin{enumerate}
\item Using degree $d$ relations reduce $A_i,B_i$ to some simple, special points of $Q$ contained in a subset $L\subset Q$. \hskip 12,2cm (*)
\item Show that any relation between the points of $L$ is generated in degree $d$.
\end{enumerate}
In general either of these two points can be very difficult.
\uwa\label{wierzcholkinie}\emph{
It is well known that the projective toric variety defined by a polytope $Q$ is covered by affine subsets given by localizations with respect to variables corresponding to vertices of $Q$. Thus one can be tempted to prove that $I_Q=I_{Qd}$ only in the localizations by variables corresponding to vertices. Note however that in general, we do not know if the scheme defined by $I_{Qd}$ is also covered by open sets given by localizations by variables corresponding to vertices. Indeed, $I_{Qd}$ and $I_Q$ may be different on the set-theoretic level. For example if $\Proj I_{Qd}$ contains a point that is zero on the coordinates corresponding to vertices and nonzero on some other coordinates, then such a point will not belong to any open set corresponding to a localization with respect to vertices. However, if $\rad I_{Qd}=I_Q$, then of course it is enough to consider localizations with respect to vertices. See also \cite{Bruns}.}
\kuwa

As the only integral points of our polytopes $P$ are the vertices, the problem described in Remark \ref{wierzcholkinie} does not concern us.
We have the following equivalences for a toric ideal $I_Q$ given by a polytope $Q$.
\begin{itemize}
\item All relations between vertices of $Q$ are generated in degree $d\Leftrightarrow$
the ideal $I$ is generated in degree $d$.
\item For any point $R\in Q$ and any relation there is an integer $m$
such that after adding $mR$ to both sides of the relation, it is generated in degree $d\Leftrightarrow$
the projective scheme defined by $I_Q$ is also defined by $I_{Qd}$.
\end{itemize}
\lem\label{justone}
Suppose that $Q$ is a group-based polytope with a group $G$ acting transitively on its integral points. The projective scheme $\Proj \c[Q]$ can be represented by an ideal generated in degree at most $d$ if and only if \emph{there exists} a point $R\in Q$ such that for any relation $\sum A_i=\sum B_i$ between points of $Q$ for $m$ sufficiently large $\sum A_i+mR=\sum B_i+mR$ is generated in degree $d$.
\klem
\dow
As $G$ acts transitively the assumption that there exists a point $R$ with the given property is equivalent to the fact that the property holds for all integral points of $Q$.
\kdow

\section{Combinatorial lemma}\label{combinatoriallemma}
\dfi[{[n]}] For any $n\in \n$ let $[n]:=\{i\in\Z: 0\leq i\leq n\}$.
\kdfi
\dfi[Coloring]
A coloring of length $n+1$ is a function $f:[n]\rightarrow[g]$. The number $g$ is called the number of colors. The support of the coloring $f$ is defined as $\{k\in[n]: f(k)\neq 0\}$.
\kdfi
\dfi[Transformation]
Consider two colorings $f_1,f_2:[n]\rightarrow[g]$. Suppose that there exist two numbers $0\leq k_1,k_2\leq n$ such that $k_j$ is \emph{not} in the support of $f_j$. Moreover suppose $f_1(k_2)=f_2(k_1)$. Define $f_j'(x)=f_j(x)$ for $x\neq k_1,k_2$. Moreover, $f_1'(k_j):=f_2(k_j)$ and $f_2'(k_j):=f_1(k_j)$. We call $f_1',f_2'$ a \emph{transformation} of $f_1$ and $f_2$.
\kdfi
Transformation of colorings corresponds to exchanging the fixed color in two colorings with the $0$ color. A multiset of colorings can be transformed into another by choosing two colorings and transforming them. We generate an equivalence relation on multisets of colorings by transformations. Abusing the notation the relation is also called transformation.
\dfi[bad coloring]
Fix $g,k,N'\in\N$.
Consider a coloring $f:[kN']\rightarrow g$. We say that $f$ is \emph{bad} (with respect to $N'$) if:
\begin{enumerate}
\item its support is contained in an interval $[(t+1)N']\setminus[tN']$ for some $0\leq t\leq k-1$,
\item it is surjective.
\end{enumerate}
Note that for $g=1$ and $f\not\equiv 0$ one has to check only the first condition.
\kdfi
\lem\label{combinatorial}
Let us fix three natural numbers: $g$ (number of colors), $s$ (bound on the support) and $a\geq 2$. Fix $\epsilon >0$. There exists $N\in \n$, such that for all $n\geq N$ any collection of colorings $f_1,\dots,f_m:[n]\rightarrow [g]$ with support of cardinality at most $s$ can be transformed into a collection $f_1',\dots,f_m'$ with the following property:

 there exist $\lfloor(1-\epsilon)\frac{n}{a}\rfloor$ numbers $x<n$ divisible by $a$, such that for any $f_j'$ and any $x$ at most one of the numbers $x,x+1,\dots,x+a-1$ is in the support of $f_j'$.
\klem
Note that the statement is nontrivial, as the number $m$ of colorings is not bounded.
\dow
The proof is inductive on $s$. For $s=1$ we can take all numbers divisible by $a$ smaller than $n$, thus the statement is trivial.

Assume that the statement is true for $s$. We will prove it for $s+1$, proceeding inductively on $g$. Fix $\epsilon$.

Suppose $g=1$. By induction hypothesis consider $N'>>0$ for which the lemma is satisfied for $s'=s$, $g'=1$, $a'=a$ and $\epsilon'=\epsilon/2$. Consider $N=k(N'+1)$ for large $k$. If we have two bad colorings (with respect to $N'+1$) with supports in different intervals $I_t$ we can transform them to two colorings that are not bad. Thus we can assume that all bad colorings have support contained in one interval, say $I_0$. By the induction hypothesis, we can transform the colorings, so that in each other interval $I_t$, we can choose $\lfloor(1-\epsilon/2)\frac{N'+1}{a}\rfloor$ numbers divisible by $a$ satisfying the required property. We see that for $k$ large enough the lemma holds.

The case for greater $g$ is only slightly more complicated. By induction, we choose $N'>>0$ for which the lemma holds for $\epsilon'=\epsilon/(g+2)$, both in case when $g$ \emph{or} $s$ is smaller. Consider $N=k(N'+1)$.
As before we can assume that all bad colorings (with respect to $N'+1$) have support contained in one interval $I_t$, say $I_0$. For any other interval $I_t$, all colorings that are not bad can be divided into (not necessarily disjoint) subsets:
\begin{itemize}
\item $g$ subsets depending on the color that is missing in the image,
\item one subset containing colorings with supports not contained in the given interval.
\end{itemize}
For each coloring that is not bad we choose one of the above $g+1$ subsets to which the coloring belongs and fix this choice. By induction we can transform the colorings in the following way: for each interval apart from $I_0$, for each subset of colorings, there are at most $\lfloor(\epsilon/(g+2))\frac{N'+1}{a}\rfloor+1$ numbers divisible by $a$ that do not satisfy the statement of the lemma. Hence, the lemma holds for $k$ large enough.
\kdow
\section{Bounding the degree of phylogenetic invariants}\label{proof2}
\dfi[Support of a group-based flow] Let $f$ be any group-based flow. The set of edges to which $f$ associates a nonneutral element is called the
\emph{support} of $f$. \kdfi
As a warmup, in the next two lemmas, we prove Main Theorem 2 for group-based models.
\lem\label{malysupport}
Consider a claw tree $K_{1,n}$ and a finite abelian group $G$. Let $\nt$ be the neutral element of $\flo$ and, abusing the notation, also the corresponding vertex of the associated polytope $P$.
Consider any relation $\sum A_i=\sum B_i$ of integral points of $P$. There exists $m$ such that the relation $m\nt+\sum A_i=m\nt+\sum B_i$ can be transformed, using only quadrics, to a relation among group-based flows with support of cardinality at most $D(G)$, where $D(G)$ is the Davenprot's constant of the group $G$. In particular, when $G$ is not a cyclic group then the supports are of cardinality at most $|G|-1$.
\klem
\dow
Consider any group-based flow $A$. Suppose its support is of cardinality greater than $D(G)$. Then we can find a proper subset $S$ of edges in the support such that $\sum_{e\in S} A(e)$ is the neutral element of $G$, where the addition is taken in $G$. Thus $\nt+A$ equals the sum of two group-based flows with strictly smaller support. The lemma follows easily.
\kdow
\lem\label{abelian}
For any finite, abelian group $G$ there exists $N$, such that if the ideal $I(K_{1,N},G)$ is generated in degree $d$, then the scheme $X(K_{1,m},G)$ can be defined by an ideal generated in degree $d$ for any $m\geq N$.
\klem
\dow
Let us choose $N$ from Lemma \ref{combinatorial} for $g=|G|-1$, $s=|G|$, $\epsilon<1/2$ and $a=2$. First let us prove the following claim:

\emph{Any relation between group-based flows of length $m$ with support of cardinality at most $|G|$ is generated in degree $d$.}

Let us prove it inductively on $m$, starting from $m=N$, where the assumption holds. Note that each group-based flow can be considered as a coloring; we number group elements from $0$ to $|G|-1$, assigning zero to the neutral element.
Transformation of two colorings corresponds to a quadric.
By Lemma \ref{combinatorial}, as $\epsilon <1/2$ we may transform the colorings to such that there exists an even number $x$, such that $x$ or $x+1$ is not in the support for any group-based flow appearing in the relation. Recall that $e_0$ is the neutral element of the group $G$.

Let us replace each group-based flow $f$ appearing in the relation by:
$$f'(j)=
\begin{cases}
f(j)\quad j<x\\
f(j+1)\quad j>x
\end{cases}
f'(x)=
\begin{cases}
f(x) \text{ if }f(x)\neq e_0\\
f(x+1)\text{ if }f(x+1)\neq e_0\\
e_0\text{ if }f(x)=f(x+1)=e_0.\\
\end{cases}
$$
We obtain a relation between group-based flows $f'$ of length $m-1$. Thus, by induction, it is generated in degree $d$. Note that each generation step can be lifted to the relation between $f$. Finally, using quadrics, we can rearrange entries indexed by $x$ and $x+1$, which proves the claim.

Due to Lemma \ref{malysupport}, by adding to any relation the group-based flow $\nt$ sufficiently many times, one can reduce it to a relation between group-based flows with support of cardinality at most $|G|$ using quadrics. As the group $\flo$ acts transitively on $P$, by Lemma \ref{justone} it is enough to consider relations multiplied by a high power of $\nt$. Thus the lemma follows.
\kdow
\uwa
The constant $N$ can be explicitly computed by explicit calculation of $k$ appearing in the proof Lemma \ref{combinatorial}. Moreover, due to standard theorems, a toric variety in bounded dimension and bounded exponents of defining monomials has got an explicitly bounded degree of generation - \cite[Theorem 4.7, Corollary 4.15]{Stks}.
\kuwa
The case of $G$-models is based on the same ideas, but is technically more involved.
\dow[Proof of Main Theorem 2]
Let us fix a tree $T=K_{1,m}$. We prove that for $m$ sufficiently large the degree in which the corresponding ideal is generated is the same as for the tree $K_{1,m-|G|}$.

Consider any relation $\sum A_i=\sum B_i$ between integral points of $P_{sub}$.

Step 1 - reducing points $A_i,B_i$ to special points.

Fix any integral point $R\in P_{sub}$. It can be lifted to a group-based flow $(g_1,\dots,g_m)$. By reordering coordinates we may assume that $g_1=\dots=g_k$ for $k=\lceil \frac{n}{|G|}\rceil$. Fix a number $i$. Consider any lift $(a^i_1,\dots,a^i_m)$ of $A_i$. Suppose there are more than $|G|^2$ elements $a^i_j\neq g_1$ for $j\leq k$. We can choose a subset $J\subset\{1,\dots,k\}$ such that:
\begin{enumerate}
\item $|J|$ is divisible by $|G|$,
\item $\sum_{j\in J} a^i_j=e_0$, where $e_0$ is the neutral element of the group $G$,
\item $a^i_j\neq g_1$ for $j\in J$.
\end{enumerate}
 We can apply a quadratic relation on $R+A_i$ replacing $a^i_j$ for $j\in J$ with corresponding $g_1$. This decreases the number of elements $a^i_j$ different from $g_1$.
Note that any relation that holds between the vertices of $P$ can be projected to a relation that holds between corresponding vertices of $P_{sub}$. Thus using quadratic relations we can assume that for each fixed lift of $B_i$ and $A_i$ on the first $k$ coordinates there are at most $|G|^2$ elements different from $g_1$.
This finishes step 1 and allows us to consider only relations between such special points.

Step 2 - generating the relation.

In the case of group-based model we could have assumed that $g_1=e_0$ due to the action of the group $\flo$. Now, we cannot make this assumption, as the action may be not compatible with the projection $P\rightarrow P_{sub}$. We will apply Lemma \ref{combinatorial} for $a=|G|+1$. Recall that we have chosen lifts of elements $A_i$ and $B_i$. We will associate colorings to each $A_i$, $B_i$ depending on the lift. Namely, the lift $g_1$ corresponds to color $0$ and other group elements to remaining $|G|-1$ colors. We restrict the domain of all the colorings to the interval $[k-1]$, as otherwise we would not be able to bound the support. As in the proof of Lemma \ref{abelian} we may find a number $0\leq x\leq k-1$, such that if we consider the lift $(l_1,\dots,l_n)$ of any $A_i$ or $B_i$ there is at most one $l_i\neq g_1$, for $i=x,x+1,\dots,x+|G|$. We will now contract elements $l_x,\dots,l_{x+|G|}$ deleting $|G|$ times the element $g_1$. Note that we will obtain a group-based flow on the tree $K_{1,m-|G|}$, as the order of $g_1$ divides $|G|$. Let $\tilde A_i$, $\tilde B_i$ be the group-based flows obtained in this way from the lifts of $A_i$ and $B_i$.
Let $\pi: M_E\rightarrow M_{E,sub}$ be the projection. As $\sum A_i=\sum B_i$ we have $\sum \pi(\tilde A_i)=\sum \pi(\tilde B_i)$, as passing from the lift of $A_i$ and $B_i$ to $\tilde A_i$ and $\tilde B_i$ we have deleted elements that projected to the same element. Suppose the relation $\sum \pi(\tilde A_i)=\sum \pi(\tilde B_i)$ is generated in degree $d$. Note\footnote{In fact this is the crucial property of a submodel of a group-based model that is necessary for the proof to work.}, that the basis elements of $M_e$ project to basis elements of $M_{e,sub}$. In particular, each coordinate $b_{e,[g]}^*$ for a fixed $e$ is the same for both sides of any relation. This allows us to lift the generation process of the relation $\sum \pi(\tilde A_i)=\sum \pi(\tilde B_i)$ to $\sum A_i=\sum B_i$. Finally, we can generate the relation using quadratic equations, exchanging the parts indexed by $x,\dots,x+|G|$.
\kdow

\section{Proof of the main theorem \ref{glKimura}}\label{proof}
The whole Section is devoted to the proof of Main Theorem \ref{glKimura}. Below we present an equivalent statement.
\tw\label{gltw}
For any tree $T$ the ideal $I(T,\z_2\times\z_2)$ and the subideal generated in degree at most four define the same projective scheme.
\ktw
Following the arguments of \cite[Chapter 5]{SS} one immediately reduces to the case when $T=K_{1,l}$ for certain $l\in \n$.
Let us restate the general definitions already introduced, in case of the group $G=\z_2\times\z_2$ (the 3-Kimura model) on a claw tree $K_{1,l}$.
\dfi[Group-based flow, ideals $I$ and $I_4$, polytope $P$]
A group-based flow is an association of elements of the group $\z_2\times\z_2$ to edges of $K_{1,l}$, such that the sum of all elements is the natural element. Let $I$ be the ideal of the variety $X(K_{1,l},\z_2\times\z_2)$ and let $I_4$ be the subideal generated in degree $4$. Let $P$ be the corresponding polytope.
\kdfi
We will identify a group-based flow with an $l$-tuple of group elements summing up to zero. The sum of such
$l$-tuples will be a coordinatewise sum, where each entry is treated as an element of the free abelian group generated
by elements of $\z_2\times \z_2$. Each group-based flow represents a vertex of a polytope $P$. The addition described above is the addition in the lattice generated by vertices of the polytope.
\ex \emph{For $l=4$ we can add:}
$$((0,0)+(0,1),(1,0)+(1,1),2(0,1),2(0,0))+((0,1),(1,0),(1,1),(0,0))$$$$=((0,0)+2(0,1),2(1,0)+(1,1),2(0,1)+(1,1),3(0,0)).$$
\emph{The neutral group-based flow} $\nt=((0,0),(0,0),(0,0),(0,0))$. \kex

\dfi[Pair, triple]
We say that a group-based flow is a \emph{pair} if and only if the cardinality of the support is equal to two.
We say that a group-based flow is a \emph{triple} if and only if the cardinality of the support is equal to three.
\kdfi
 By Lemmas \ref{justone} and \ref{malysupport} we have to generate relations only between group-based flows that are pairs and triples. This completes the first step of the method (*) presented in Section \ref{idea}. The set $L$ consists of pairs and triples. Note that this part of the proof can be adjusted to other groups $G$, for $L$ consisting of group-based flows with support of cardinality at most $|G|$.

Let us fix any relation $\sum A_i=\sum B_i$, where $A_i$ and $B_i$ are either pairs or triples.
Our aim is to transform $\sum A_i$ to $\sum B_i$ in a series of steps, each time replacing at most four $A_i$ by group-based flows with the same sum\footnote{We are also allowed to add the group-based flow $\nt$ to both sides.}. We assume that among $A_i$ there are more or the same number of triples as among $B_i$. We first
try to reduce the relation, so that consequently:
\begin{enumerate}
\item Among $A_i$ there are as few triples as possible,
\item Among $B_i$ there are as few triples as possible,
\item The degree of the relation is as small as possible.
\end{enumerate}
More precisely let $t$ and $t'$ be the number of triples among respectively $A_i$ and $B_i$. Let $d'$ be the degree of the relation.
Our proof will be inductive on $(t,t',d')$ with lexicographic order.

To prove Theorem \ref{gltw} we consider separately three cases depending on the number of triples among $A_i$. The cases
are:

a) there are no triples,

b) there is exactly one triple,

c) there are at least two triples.


We say that a family of group-based flows agrees on an index $j$ of an edge if they all associate the same element to $j$ and $j$ belongs
to their \emph{support}. We will denote by $g_1$, $g_2$ and $g_3$ the three nonneutral elements of $\z_2\times\z_2$. A triple that
associates $g_1$ to index $a$, $g_2$ to index $b$ and $g_3$ to index $c$ is denoted by $(a,b,c)$. A pair that associates an element $g_i$ to indices
$d$ and $e$ will be denoted by $(d,e)_{g_i}$ and called a $g_i$ pair. We say that $g_i$ is contained in a group-based flow if there exists an index $j$, such that the group-based flow associates $g_i$ to $j$.
We believe that the following proofs are impossible to follow without a piece of paper. We strongly encourage the reader to note what group-based flows appear on both sides of the relation at each step of the proof.
\subsection{The case with no triples}
First note that there are no triples among $B_i$. Without loss of generality we may assume that $A_1$ is a pair equal to $(a,b)_{g_1}$. 
Hence there exists $(b,c)_{g_1}$ among $B_i$ for some index $c$. If $c=a$ we can reduce
this pair, hence we assume $c\neq a$. There exists a group-based flow, say $A_2$ that is $(c,d)_{g_1}$. If $d=b$ we can reduce this
pair. We consider two other cases:

1) $d\neq a$. Then we use the degree two relation $(a,b)_{g_1}+(c,d)_{g_1}=(a,d)_{g_1}+(b,c)_{g_1}$ and we can reduce $(b,c)_{g_1}$.

2) $d=a$. Then there is a group-based flow, say $B_1$ given by $(a,e)_{g_1}$. If $e=b$ or $e=c$ we can reduce this pair. In the other
cases we use the relation $(a,e)_{g_1}+(b,c)_{g_1}=(a,b)_{g_1}+(e,c)_{g_1}$ and we reduce $(a,b)_{g_1}$.

Notice that in this very easy case we have only used degree two relations.
\subsection{The case with one triple}
Let $A_1$ be the only triple among $A_i$.
\lem
There is exactly one triple among $B_i$.
\klem
\dow
Due to the assumptions we know that there is at most one triple among $B_i$. We exclude the case when there are no triples by comparing the parity of the number of times the element $g_1$ appears on both sides of the relation.
\kdow
Due to the previous lemma we may assume that $B_1$ is the only triple among $B_i$.
Without loss of generality, assume $A_1=(1,2,3)$.

\subsubsection{Case: the triples agree on at least two elements in their support}\label{dwadwa}
Suppose that $A_1=(1,2,3)$ and $B_1=(1,2,c)$. Of course if $c=3$ we can make a reduction. In other case we must have a pair
$(c,d)_{g_3}$ among $A_i$. If $d\neq 3$ then we use the relation $(c,d)_{g_3}+(1,2,3)=(1,2,c)+(3,d)_{g_3}$ and reduce the triples.
Assume $d=3$. Analogously, we can assume there is a pair $(3,c)_{g_3}$ among $B_i$, hence we can reduce this pair.
\subsubsection{Case: the triples agree on exactly one element in their support}\label{oneone}
Consider the case when triples agree on at least one element, say $1$, in their common support. By the previous case we may assume that they agree on exactly one element.

As before
let $A_1=(1,2,3)$ and $B_1=(1,b,c)$. We consider three cases.

1) $b\neq 3$.

There must be a pair $(b,d)_{g_2}$ among $A_i$. If $d\neq 2$ then we can apply the relation
$(b,d)_{g_2}+A_1=(1,b,3)+(d,2)_{g_2}$. This reduces to the case \ref{dwadwa}. So we assume
$d=2$. There must be a pair $(2,e)_{g_2}$ among $B_i$. We may assume $e\neq b$ as otherwise we would be able to make a reduction. 
Hence there
must also be a pair $(e,f)_{g_2}$ among $A_i$. If $f\neq b$ we can use a relation $(e,f)_{g_2}+(2,b)_{g_2}=(e,2)_{g_2}+(f,b)_{g_2}$ and reduce
$(e,2)_{g_2}$. For $f=b$ we must have a pair $(b,g)_{g_2}$ among $B_i$. If $g=2$ or $g=e$ then this pair can be reduced.
In the other case we use the relation $(e,2)_{g_2}+(b,g)_{g_2}=(e,g)_{g_2}+(b,2)_{g_2}$ and reduce $(b,2)_{g_2}$.

2) $c\neq 2$.

This case is analogous to 1).

3) $b=3$ and $c=2$.

\lem\label{parazahacza}
If there is a pair $(p,q)_{g_2}$ among $A_i$, such that $p,q\neq 2$ then we may assume that it is equal to $(1,3)$.
\klem
\dow
Suppose that
$p\neq 1,2,3$ and $q\neq 2$. We apply a relation $(p,q)_{g_2}+A_1=(1,p,3)+(q,2)_{g_2}$ and reduce to case 2) $c\neq 2$.
\kdow
Analogously if there is a pair
$(p,q)_{g_2}$ among $B_i$, such that $p,q\neq 3$ then we can assume it is equal to $(1,2)_{g_2}$.

Notice that there must be a pair
$(3,d)_{g_2}$ among $A_i$ and a pair $(2,e)_{g_2}$ among $B_i$. From Lemma \ref{parazahacza}, $d$ equals either $2$ or $1$
and $e$ equals either $3$ or $1$. We will consider subcases.

3.1) Suppose that $d=2$.

If $e=3$ then we can make a reduction of pairs. If $e=1$ we must have a pair $(1,f)_{g_2}$ among $A_i$.
If $f=2$ we make a reduction, hence we assume $f=3$. This means that there must be a pair $(3,g)_{g_2}$ among $B_i$.
If $g=2$ or $g=1$ we can make a reduction. Otherwise we apply the relation $(1,2)_{g_2}+(3,g)_{g_2}=(1,3)_{g_2}+(2,g)_{g_2}$ and reduce the pair $(1,3)_{g_2}$.

3.2) Suppose that $e=3$.

This case is similar to $3.1)$.

3.3) Suppose that $d=1$ and $e=1$.

As this is the only case left we may repeat the same reasoning for $g_3$. In particular, we can assume there is a pair $(1,2)_{g_3}$ among $A_i$.
We see that we can reduce the triples by applying the
following relation: $$(1,2,3)+(1,3)_{g_2}+(1,2)_{g_3}=(1,3,2)+(1,2)_{g_2}+(1,3)_{g_3}.$$ This is a degree three
relation.
\subsubsection{Case: The triples do not agree on any element of the support}\label{triplesdonotagree}
We want to reduce to one of previous cases. 
Consider the following two cases.

1) \emph{The triples $A_1$ and $B_1$ have different supports.}

Once again let $(1,2,3)=A_1$ and let $(a,b,c)=B_1$. We may assume that $a$ is not in the support of $A_1$. We see that there must
be a pair $(a,f)_{g_1}$ among $A_i$. If $f\neq 1$ we can use a relation $(a,f)_{g_1}+A_1=(a,2,3)+(f,1)_{g_1}$. This
reduces to the case \ref{oneone}, hence we assume that $f=1$. There must be a pair $(g,1)_{g_1}$
among $B_i$. If $g=a$ we can reduce this pair, so we assume $g\neq a$. 
Notice that there must be a pair $(g,h)_{g_1}$
among $A_i$. If $h\neq a,$ then we can use relation $(1,a)_{g_1}+(g,h)_{g_1}=(g,1)_{g_1}+(h,a)_{g_1}$ and reduce the pair $(g,1)_{g_1}$. So we can
assume $h=a$. Then there must be a pair $(a,i)_{g_1}$ among $B_i$. If $i=1$ then we can reduce it. Otherwise we can use the relation $(g,1)_{g_1}+(a,i)_{g_1}=(g,a)_{g_1}+(1,i)_{g_1}$ and reduce the pair $(g,a)_{g_1}$.

2) \emph{The set $\{1,2,3\}$ is the support of $B_1$ and $A_1$.}

Remember that due to the assumption \ref{triplesdonotagree} the triples $A_1$ and $B_1$ do not agree on any element from their support. Without loss of generality we may assume $A_1=(1,2,3)$ and $B_1=(2,3,1)$. Hence there must be a pair $(2,a)_{g_1}$ among $A_i$ and $(1,b)_{g_1}$ among $B_i$. If $a=1$ and $b=2$ then both pairs are the same and can be reduced. As both cases are symmetric we can assume that $a\neq 1$.

If $a\neq 3$ we can use the relation $(2,a)_{g_1}+(1,2,3)=(a,2,3)+(2,1)_{g_1}$. This reduces to the case with
different supports. We are left with the case $a=3$. There must be a pair $(3,z)_{g_1}$ among $B_i$. If $z\neq 1$ we can use the
relation $(3,z)_{g_1}+B_1=(z,3,1)+(2,3)_{g_1}$. This would enable to reduce the $(2,3)_{g_1}$ pair and decrease the degree.
So we can assume that $z=1$. So far we have shown that there must be pairs $(2,3)_{g_1}$ among $A_i$ and $(3,1)_{g_1}$ among $B_i$\footnote{Notice that we have made a symmetry assumption $a\neq 1$. The symmetric assumption would be $b\neq 2$. However as the result we got was symmetric, also for $b\neq 2$ we prove the existence of the same pairs.}.
By the same reasoning for $g_2$ and $g_3$ we see that we can use the following relation:
$$(1,2,3)+(2,3)_{g_1}+(1,3)_{g_2}+(1,2)_{g_3}=(2,3,1)+(2,3)_{g_3}+(1,3)_{g_1}+(1,2)_{g_2}.$$
Notice that this is a degree four relation. It enables us to reduce triples.

\subsection{The case with at least two triples}
We suppose that there are at least two triples among $A_i$.
\lem\label{trojkizgadzajasienajakims}
If there are two triples $A_1$, $A_2$ among $A_i$ that do not agree on any element of their supports then we can make a reduction. Thus we can assume that any two triples among $A_i$ agree on at least one index.
\klem
\dow
The assumptions are equivalent to $A_1=(a,b,c)$, $A_2=(d,e,f)$ with $a\neq d$, $b\neq e$, $c\neq f$. We apply the relation $A_1+A_2+\nt=(a,d)_{g_1}+(b,e)_{g_2}+(c,f)_{g_3}$ that reduces the number of triples.
\kdow
\lem\label{notallagree}
If there is no index on which all triples among $A_i$ agree then we can make a reduction.
\klem
\dow
Suppose there is no index on which all triples among $A_i$ agree. We may consider only two cases due to Lemma \ref{trojkizgadzajasienajakims}.

1) Suppose that any two triples from $A_i$ agree on at least two elements.

Consider any triple $A_1=(1,2,3)$. Due to the fact that not all triples from $A_i$ associate $g_1$ to $1$
there is a triple $(a,2,3)$ with $a\neq 1$ among $A_i$. There also must be a triple that does not associate
$g_2$ to $2$. But this cannot happen as the triple must agree with both $(1,2,3)$ and $(a,2,3)$ on two indices, which gives a contradiction.

2) There exist two triples that agree only on one index.

Let $A_1=(1,2,3)$ and $A_2=(1,b,c)$ with $b\neq 2$ and $c\neq 3$. Due to the case assumption
there is a triple $A_3=(d,e,f)$ with $d\neq 1$. Remember that any two triples have to agree on at least one element due to Lemma \ref{trojkizgadzajasienajakims}. Hence without loss of generality we can assume $e=b$ and $f=3$. We can apply the relation:
$$A_1+A_2+A_3+\nt=(d,1)_{g_1}+(2,b)_{g_2}+(3,c)_{g_3}+(1,b,3).$$
This relation reduces the number of triples.
\kdow


\emph{Due to the previous lemma we may assume that there exists an index, say $1$, such that all triples among $A_i$ associate to it the same nonneutral element, say $g_1$.}
\dfi[$k$]
Let $k$ be the number of indices on which all triples among $A_i$ agree. We know that $1\leq k \leq 3$.
\kdfi
We proceed inductively on $k$, as for $k=0$ we already know from Lemma \ref{notallagree} how to reduce the relation. Hence from now on decreasing $k$ is also a reduction.
\lem\label{postacdwojek}
Suppose that all triples $A_i$ associate $g_j$ to an index $l$. If there is a pair $(x,y)_{g_j}$ among $A_i$ with $l\neq x,y$ then either $\{l,x,y\}$ is the support of all triples among $A_i$ or we can make a reduction.
\klem
\dow
Suppose that there is a triple $A_i$ with the support $\{l,b,c\}$ different from $\{l,x,y\}$. We can assume $x\neq b,c$. We apply the relation $A_i+(x,y)_{g_j}=\tilde A_i+(l,y)_{g_j}$, where $\tilde A_i$ associates $g_j$ to $x$ and agrees with $A_i$ on $b$ and $c$. This relation reduces $k$.
\kdow
\lem\label{niewszystkozahacza}
Suppose that all triples from $A_i$ associate $g_j$ to an index $l$. If all pairs $(x,y)_{g_j}$ among $A_i$ have $l$ in the support then we can reduce all such pairs.
\klem
\dow
Recall that $t$ is the number of triples among $A_i$. Let $p$ be the number of $g_j$ pairs among $A_i$. Let $t'_1$ (resp. $t'_2$) be the number of triples in $B_i$ that assign (resp. do not assign) $g_j$ to $l$. Let $p'_1$ (resp. $p'_2$) be the number of $g_j$ pairs among $B_i$ that have (resp. do not have) $l$ in the support. We know that $t\geq t'_1+t'_2$. Comparing the number of times $g_j$ appears in $A_i$ and $B_i$ we get:
$$t+2p=t'_1+t'_2+2(p_1'+p_2').$$
Comparing the number of times $g_j$ appears on index $l$ we get:
$$t+p=t'_1+p_1'.$$
This forces $t'_2=p_2'=0$, $t=t'_1$ and $p=p'_1$. Hence all $g_j$ pairs and triples among $A_i$ and $B_i$ must assign $g_j$ to $l$. Hence the multisets of pairs must be the same for $A_i$ and $B_i$.
\kdow
\lem
Suppose that all triples from $A_i$ associate $g_j$ to an index $l$.
If there are $g_l$ pairs among $A_i$, then we can make a reduction.
\klem
\dow
Without loss of generality we assume $g_l=g_1$.
Due to Lemma \ref{niewszystkozahacza}, it is enough to prove that if there are pairs $(a,b)_{g_1}$ among $A_i$ with $a,b\neq 1$ then we can make a reduction. Suppose that there is such a pair. 
Due to Lemma \ref{postacdwojek} all the triples among $A_i$ must have the support $\{1,a,b\}$. So either $k=1$ or $k=3$. If $k=1$ we can apply the relation $$(1,a,b)+(1,b,a)+(a,b)_{g_1}+\nt=(1,a)_{g_1}+(1,b)_{g_1}+(a,b)_{g_2}+(a,b)_{g_3}.$$ This reduces the number of triples.
Thus we can assume that all triples among $A_i$ are equal to $(1,a,b)$.

Claim:\emph{
Consider any pair $(c,d)_{g_2}$ among $A_i$.
We can assume that its support is contained in $\{1,a,b\}$.}

\dow[Proof of the Claim]
 Suppose this is not the case, that is $c\not\in \{1,a,b\}$. Due to Lemma \ref{postacdwojek} we can assume $d=a$.

 1) Suppose that there is a $g_2$ pair among $A_i$ that does not contain $a$ in the support.

  It must be equal to $(1,b)_{g_2}$ due to Lemma \ref{postacdwojek}. We can apply the relation $(1,b)_{g_2}+(a,c)_{g_2}=(c,1)_{g_2}+(a,b)_{g_2}$. Applying once again Lemma \ref{postacdwojek} to the pair $(c,1)_{g_2}$ we can make a reduction.

2) All $g_2$ pairs among $A_i$ contain $a$ in the support.

Due to Lemma \ref{niewszystkozahacza} we can make a reduction.
\kdow
Thus the support of all $g_2$ pairs among $A_i$ is contained in $\{1,a,b\}$. The same holds for $g_1$ and $g_3$ pairs. Thus all group-based flows among $A_i$ have support contained in $\{1,a,b\}$. Hence the same must hold for $B_i$. So our relation is a relation only on three indices. It is well-known \cite{SS} that the ideal for a tree with three edges is generated in degree $4$, so in particular the considered relation is generated in degree $4$.
\kdow
\wn\label{glwn}
If all triples among $A_i$ associate $g_j$ to an index $l$, then there are no $g_j$ pairs among $A_i$. Consequently, there are no $g_j$ pairs among $B_i$ and all triples among $B_i$ associate $g_j$ to $l$. Moreover, the number of triples among $A_i$ equals the number of triples among $B_i$.\kwadrat
\kwn
By the previous corollary \emph{we assume} that there are no $g_1$ pairs neither among $A_i$ nor $B_i$. Moreover, there is the same number of triples among $A_i$ and $B_i$ and they all associate $g_1$ to $1$.
\lem\label{differentsupport}
If all the triples among $A_i$ and $B_i$ have support contained in $\{1,2,3\}$ then we can make a reduction.
\klem
\dow
Suppose all triples have support contained in $\{1,2,3\}$.
In this case $k=1$ or $k=3$. If $k=1$ then among $A_i$ there is a triple $(1,2,3)$ and $(1,3,2)$. Any triple among $B_i$ is equal to one of those. In particular, one of these triples can be reduced. If $k=3$ there are no pairs. All triples among $A_i$ and $B_i$ are equal, thus the relation is trivial.
\kdow

\subsubsection{Case: $k=1$}
We first consider the most difficult case $k=1$. As always let $A_1=(1,2,3)$ and $B_1=(1,b,c)$.
As the proof is quite complicated we decided to include the diagram that describes most important cases. While reading the proof we encourage the reader to follow at which node we are. The proof is "depth-first, left-first".

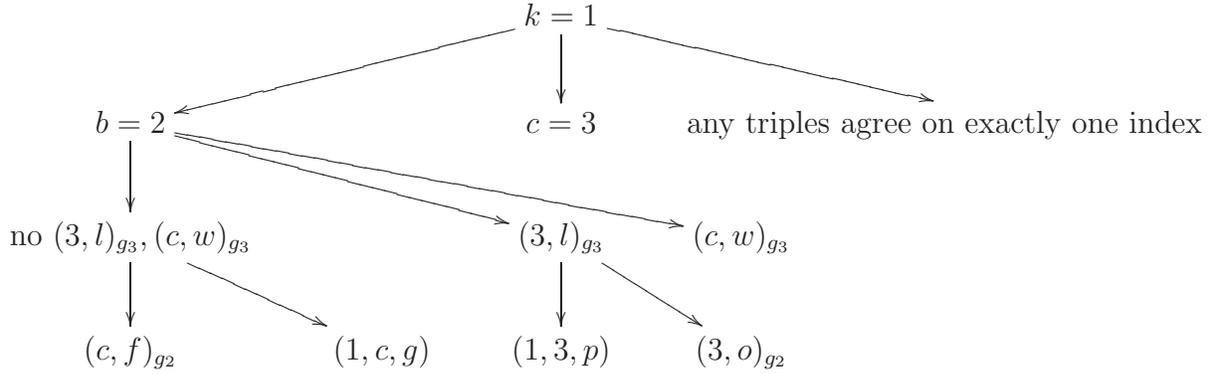
\[
\xymatrix{
&&k=1\ar[lld]\ar[d]\ar[rrd]&&\\
b=2\ar[d]\ar[drr]\ar[drrr]&&c=3&&\hskip -70pt\text{any triples agree on exactly one index}\\
\text{no }(3,l)_{g_3},(c,w)_{g_3}\ar[d]\ar[dr]&&(3,l)_{g_3}\ar[d]\ar[dr]&(c,w)_{g_3}&\\
(c,f)_{g_2}&(1,c,g)&(1,3,p)&(3,o)_{g_2}&\\
}
\]

We start with the left node in the second row -- assume $b=2$.  Then we may assume $c \neq 3$, or else we can reduce.

We move to the most left node in the third row --
 suppose that there is no $g_3$ pair among $A_i$ that has got $c$ in the support and, symmetrically, there is no $g_3$ pair among $B_i$ that has got $3$ in the support. There must be a triple $(1,e,c)$ among $A_i$. If $e\neq 3$ then we apply the relation $(1,2,3)+(1,e,c)=(1,2,c)+(1,e,3)$ and reduce the triple $(1,2,c)$. We may assume $e=3$. Analogously, we may assume that there is a triple $(1,c,3)$ among $B_i$. Hence there must be either a pair $(c,f)_{g_2}$ or a triple $(1,c,g)$ among $A_i$.

We continue to the most left node in the fourth row -- suppose that there is a pair $(c,f)_{g_2}$. If $f\neq 2$ we apply the relation $(1,2,3)+(c,f)_{g_2}=(1,c,3)+(f,2)_{g_2}$ and reduce the triple $(1,c,3)$. If $f=2$ we apply the relation $(1,3,c)+(c,2)_{g_2}=(1,2,c)+(3,c)_{g_2}$ and reduce the triple $(1,2,c)$.

Hence we can assume that there is a triple $(1,c,g)$ among $A_i$ -- second node in the fourth row. If $g\neq 2$ then we apply the relation $(1,c,g)+(1,2,3)=(1,2,g)+(1,c,3)$ and reduce the triple $(1,c,3)$. For $g=2$ we apply the relation $(1,2,3)+(1,3,c)+(1,c,2)=(1,2,c)+(1,3,2)+(1,c,3)$ and reduce the triple $(1,2,c)$.

 We continue to the second node in the third row. We assume that there is a pair $(3,l)_{g_3}$ among $B_i$. If $l\neq c$ we apply the relation $(1,2,c)+(3,l)_{g_3}=(1,2,3)+(c,l)_{g_3}$ and reduce the triple $(1,2,3)$. If there was a pair $(c,m)_{g_3}$ among $A_i$ then analogously we could assume $m=3$ and we would be able to reduce this pair. So there must be a triple $(1,n,c)$ among $A_i$. If $n\neq 3$ then we apply the relation $(1,2,3)+(1,n,c)=(1,n,3)+(1,2,c)$ and reduce the triple $(1,2,c)$. So we assume $A_2=(1,3,c)$. Hence there is either a pair $(3,o)_{g_2}$ or a triple $(1,3,p)$ among $B_i$.

 We move to the third node in the fourth row -- suppose that there is a triple $(1,3,p)$ among $B_i$. If $p\neq 2$ we apply the relation $(1,2,c)+(1,3,p)=(1,2,p)+(1,3,c)$ and we reduce $(1,3,c)$. So we assume $p=2$.
 We apply the relation $(1,3,2)+(3,c)_{g_3}=(1,3,c)+(2,3)_{g_3}$ and reduce the triple $(1,3,c)$.



We pass to the fourth node in the fourth row -- we assume that there is a pair $(3,o)_{g_2}$ and there is no triple $(1,3,p)$ among $B_i$. If $o\neq 2$ then we apply the relation $(1,2,c)+(3,o)_{g_2}=(1,3,c)+(2,o)_{g_2}$ and reduce $(1,3,c)$. So we assume there is a pair $(2,3)_{g_2}$ among $B_i$. Suppose that this pair appears $r>0$ times among $B_i$. Note that we may assume that there are no pairs $(2,s)_{g_2}$ among $A_i$. Indeed suppose that there is such a pair. If $s\neq 3$ then we apply the relation $(1,3,c)+(2,s)_{g_2}=(1,2,c)+(3,s)_{g_2}$ and reduce the triple $(1,2,c)$. If $s=3$ we reduce the pair $(2,3)_{g_2}$. Hence we assume there are at least $r+1$ triples of the type $(1,2,t)$ among $A_i$. If there is a triple with $t\neq 3$ then we apply the relation $(1,3,c)+(1,2,t)=(1,3,t)+(1,2,c)$ and reduce the triple $(1,2,c)$. Hence we assume there are at least $r+1$ triples $(1,2,3)$ among $A_i$. Notice that we may assume there are no triples of the type $(1,y,3)$ among $B_i$. Indeed, in such a case we could apply the relation $(1,y,3)+(2,3)_{g_2}=(1,2,3)+(y,3)_{g_2}$ and reduce $(1,2,3)$. Hence we assume there are at least $r+1$ pairs of the type $(3,u)_{g_3}$ among $B_i$. If $u\neq c$ then we apply the relation $(1,2,c)+(3,u)_{g_3}=(1,2,3)+(c,u)_{g_3}$ and reduce the triple $(1,2,3)$. Hence we assume there are at least $r+1$ pairs $(3,c)_{g_3}$ among $B_i$. Note that we can assume there are no pairs of the type $(c,v)_{g_3}$ among $A_i$. Indeed if $v=3$ we could reduce this pair. If $v\neq 3$ then we apply the relation $(1,2,3)+(c,v)_{g_3}=(1,2,c)+(3,v)_{g_3}$ and reduce the triple $(1,2,c)$. Hence we must have at least $r+1$ triples of the type $(1,z,c)$ among $A_i$. If $z\neq 3$ then we apply the relation $(1,2,3)+(1,z,c)=(1,2,c)+(1,z,3)$ and reduce the triple $(1,2,c)$. So we may assume there are at least $r+1$ triples $(1,3,c)$ among $A_i$. Note that the elements $g_2$ on $3$ cannot be reduced -- among $B_i$ there are only $r$ pairs containing them and no triples. The contradiction finishes this case.

Consider the third node in the third row -- there is a pair $(c,w)_{g_3}$ among $A_i$. This is completely analogous to the second node in this row, which was already considered.

Also the second node in the second row -- $c=3$ -- is analogous to the first node in the second row.

We are left with the last, third node in the second column -- any two triples $A_i$ and $B_j$ agree on exactly one index, that is on $1$.
Due to Lemma \ref{differentsupport}, there is a triple among $B_i$, say $B_1$, with support different then some triple in $A_i$, say $A_1$. Exchanging $g_2$ and $g_3$ if necessary, we can assume $b\neq 2$ and $b\neq 3$. Due to the case assumption there must be a pair $(b,d)_{g_2}$ among $A_i$. If $d\neq 2$ then we apply the relation $(1,2,3)+(b,d)_{g_2}=(1,b,3)+(d,2)_{g_2}$ and reduce to the case $b=2$\footnote{Notice that we do not reduce to the case $k=2$ as if this was true we would have already been in the first node in the second column $b=2$.}. Analogously we must have the same pair among $B_i$ and it can be reduced.

\subsubsection{Case: $k=2$ or $k=3$}
Suppose now that $k=2$.
Let $A_1=(1,2,3)$ and $B_1=(1,2,c)$. If we cannot reduce $B_1$ then there must be a pair $(c,d)_{g_3}$ among $A_i$ and a pair $(3,e)_{g_3}$ among $B_i$. If $d=3$ and $e=c$ we can reduce the pairs. Thus we can assume that $d\neq 3$. We apply the relation $(1,2,3)+(c,d)_{g_3}=(1,2,c)+(3,d)_{g_3}$ and reduce the triple $(1,2,c)$.

The last, easiest case is $k=3$. Then all triples are equal to $(1,2,3)$  and there are no pairs due to Corollary \ref{glwn}. Hence we can reduce the triples. This finishes the proof of Theorem \ref{gltw}.
\section{Open problems}\label{open}
We have already recalled two conjectures of Sturmfels and Sullivant (Conjectures \ref{Kimura} and \ref{ogolna} from the Introduction). One can state simpler conjectures, that we find important.
\hip\label{hipnowageneralna}
For any group $G$ and any tree $T$ the ideal $I(X(G,T))$ and the subideal generated in degree $|G|$ define the same projective scheme.
\khip
We believe that the methods of Section \ref{proof} could be adopted to prove cases of this conjecture, especially for groups of small cardinality.

Consider a claw tree $K_{1,n}$. Subdivide the set of leaves into two separate subsets. By adding an inner edge that separates this set one obtains a tree $T$ as shown below.
$$\resizebox{150 pt}{1.5 pt}
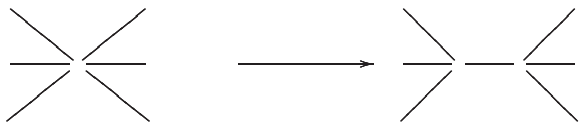
{\xymatrix{
\textbf{\ar@{-}[dr]}&&\textbf{\ar@{-}[dl]}&&&\\
&\textbf{\ar@{-}[dr]}\textbf{\ar@{-}[r]}\textbf{\ar@{-}[l]}\textbf{\ar@{-}[dl]}&&\textbf{\ar@{->}[rr]}&&\\
&&&&&\\}}
\resizebox{75 pt}{1.5 pt}{\xymatrix{
\textbf{\ar@{-}[dr]}&&&\\
&\textbf{\ar@{-}[r]}\textbf{\ar@{-}[l]}\textbf{\ar@{-}[dl]}&\textbf{\ar@{-}[r]}\textbf{\ar@{-}[dr]}\textbf{\ar@{-}[ur]} &\\
&&&\\}}
$$
One can easily see that for any\footnote{In fact for any model.} group $G$ we have $X(K_{1,n},G)\subset X(T,G)$. As long as the subsets of leaves are of cardinality at least $2$ the procedure decreases the maximal degree of the graph and allows for an inductive generation of phylogenetic invariants.
\hip[\cite{DBM} Conjecture 3.6]
The variety $X(K_{1,n},G)$ is a scheme theoretic intersection of all $X(T,G)$, where $T$ is a prolongation of $G$ with all vertices of degree at most $n-1$.
\khip
This conjecture would imply the following.
\hip
The function $\phi(G,\cdot)$ is constant.
\khip
In view of the previous conjecture and Conjecture \ref{ogolna} it is natural to bound $\phi(G,K_{1,3})$. It is an open problem if $\phi(G,K_{1,3})\leq |G|$. On the other hand, by a general theorem of Sturmfels \cite[Theorem 4.7]{Stks} we have an exponential bound on $\phi(G,K_{1,3})$. The following question, as far as we know, is open and would be an interesting step towards conjectures of Sturmfels and Sullivant.
\hip
There exists a polynomial $Q$ such that $\phi(G,K_{1,3})\leq Q(|G|)$.
\khip
One obtains a much more difficult open problem by replacing $3$ in the previous conjecture by any (fixed or even non fixed) number $n\in \N$.

The following conjecture concerns normality of the variety $X(T,G)$.
\hip
The variety $X_\p(T,\z_2\times\z_2)$ is projectively normal.
\khip
The previous conjecture is known for trees of small maximal degree. Note that $X(K_{1,3},\z_6)$ is not normal \cite[Computation 4.3]{DBM}. However, analogous questions for many other groups are open. For the group $G=\z_2$ the reader can consult \cite{DBM}.

For group-based models the cardinality of the group $G$ equals the dimension of the vector space $W$ associated to any node. For the general Markov model, the dimension of the space $W$ associated to the node, determines which secant of the Segre variety we consider. It is well--known (and easy to prove) that there are no equations of degree less or equal to $k$ for the $k$-th secant.

A positive answer to the first part of the following question would be an ideal theoretic version of the main result in \cite{DK2}. The second part of the question is in the flavor of Conjecture \ref{hipnowageneralna}, however we cannot speculate about the bound -- we only know that it is not equal to $k$.
\begin{question}
For any $k$, does there exist $d$ such that the ideal of the $k$-th secant variety of any Segre variety is generated in degree at most $d$? If yes, what is the bound on $d$ with respect to $k$?
\end{question}
The affine cone over the Segre variety is naturally embedded in a tensor product $W_1\otimes\dots\otimes W_n$. A tensor product of $n$ vector spaces can be regarded as a tensor product of $n'<n$ vector spaces by considering a tensor product of a few vector spaces as one vector space. For example:
$$W_1\otimes W_2\otimes W_3\otimes W_4=(W_1\otimes W_2)\otimes (W_3\otimes W_4).$$
Such a procedure is called a flattening. If we now the equations for the $k$-th secant variety for some $n'$, then we may induce some equations for the $k$-th secant variety for any $n>n'$, by flattenings. In the following question, it would be tempting to conjecture a linear bound.
\begin{question}
For any fixed $k$, does there exist $n'$ such that for any $n>n'$ the equations obtained by flattenings generate the ideal of the $k$-th secant variety of the Segre embedding of the product of $n$ spaces?
\end{question}

\bibliographystyle{amsalpha}
\bibliography{xbibbb}

\def\cprime{$'$}
\providecommand{\bysame}{\leavevmode\hbox to3em{\hrulefill}\thinspace}
\providecommand{\MR}{\relax\ifhmode\unskip\space\fi MR }
\providecommand{\MRhref}[2]{%
  \href{http://www.ams.org/mathscinet-getitem?mr=#1}{#2}
}
\providecommand{\href}[2]{#2}
\begin{thebibliography}{BBKM13}

\bibitem[AR08]{AllRhMarkov}
Elizabeth Allman and John Rhodes, \emph{Phylogenetic ideals and varieties for
  the general {Markov} model}, Advances in Applied Mathematics \textbf{40(2)}
  (2008), 127--148.

\bibitem[BBKM13]{BBKM}
Weronika Buczy{\'n}ska, Jaros{\l}aw Buczy{\'n}ski, Kaie Kubjas, and Mateusz
  Michalek, \emph{Degrees of generators of phylogenetic semigroups on graphs},
  arXiv:1004.1183, to appear in Central European Journal of Mathematics (2013).

\bibitem[Bru11]{Bruns}
Winfried Bruns, \emph{The quest for counterexamples in toric geometry},
  arXiv:1110.1840 (2011).

\bibitem[BW07]{BW}
Weronika Buczy\'nska and Jaros{\l}aw Wi\'{s}niewski, \emph{On geometry of
  binary symmetric models of phylogenetic trees}, J. Eur. Math. Soc.
  \textbf{9(3)} (2007), 609--635.

\bibitem[CP07]{Sonja}
Julia Chifman and Sonja Petrovi\'{c}, \emph{Toric ideals of phylogenetic
  invariants for the general group-based model on claw trees $k_{1,n}$},
  Proceedings of the 2nd international conference on {Algebraic} biology
  (2007), 307--321.

\bibitem[DBM]{DBM}
Maria Donten-Bury and Mateusz Micha{\l}ek, \emph{Phylogenetic invariants for
  group-based models}, Journal of Algebraic Statistics \textbf{3}, no.~1.

\bibitem[DE]{DE}
Jan Draisma and Rob Eggermont, \emph{Finiteness results for abelian tree
  models}, arXiv:1207.1282v1.

\bibitem[DK09]{DK}
Jan Draisma and Jochen Kuttler, \emph{On the ideals of equivariant tree
  models}, Mathematische Annalen \textbf{344(3)} (2009), 619--644.

\bibitem[DK11]{DK2}
\bysame, \emph{Bounded-rank tensors are defined in bounded degree},
  arXiv:1103.5336v2 (2011).

\bibitem[Dra10]{Dfinit}
Jan Draisma, \emph{Finiteness for the {$k$}-factor model and chirality
  varieties}, Adv. Math. \textbf{223} (2010), no.~1, 243--256.

\bibitem[ERSS04]{4aut}
Nicholas Eriksson, Kristian Ranestad, Bernd Sturmfels, and Seth Sullivant,
  \emph{Phylogenetic algebraic geometry}, Projective Varieties with Unexpected
  Properties; Siena, Italy (2004), 237--256.

\bibitem[Ful93]{Fult}
William Fulton, \emph{Introduction to toric varieties}, Annals of Mathematics
  Studies, vol. 131, Princeton University Press, Princeton, NJ, 1993, The
  William H. Roever Lectures in Geometry.

\bibitem[HMdC13]{Abraham}
Christopher Hillar and Abraham Mart{\'{\i}}n~del Campo, \emph{Finiteness
  theorems and algorithms for permutation invariant chains of {L}aurent lattice
  ideals}, J. Symbolic Comput. \textbf{50} (2013), 314--334.

\bibitem[HP89]{HendyPenny}
Michael Hendy and David Penny, \emph{A framework for the quantitative study of
  evolutionary trees}, Systematic Zoology \textbf{38} (1989), 297--309.

\bibitem[HS12]{Sullfinite}
Christopher Hillar and Seth Sullivant, \emph{Finite gröbner bases in infinite
  dimensional polynomial rings and applications}, Advances in Mathematics
  \textbf{229} (2012), no.~1, 1 -- 25.

\bibitem[JC69]{Jukes.Cantor1969}
Thomas Jukes and Charles Cantor, \emph{{Evolution of protein molecules}},
  Mammalian protein metabolism, vol. III, Academic Press, 1969, pp.~21--132.

\bibitem[Kim80]{2Kimura}
Motoo Kimura, \emph{A simple method for estimating evolutionary rates of base
  substitutions through comparative studies of nucleotide sequences}, Journal
  of Molecular Evolution \textbf{16} (1980), 111--120.

\bibitem[Kim81]{3Kimura}
\bysame, \emph{Estimation of evolutionary distances between homologous
  nucleotide sequences}, Proceedings of the National Academy of Sciences
  \textbf{78} (1981), no.~1, 454--458.

\bibitem[Man09]{Man}
Christopher Manon, \emph{The algebra of {Conformal} {Blocks}}, arXiv:0910.0577
  (2009).

\bibitem[Man12]{Manonneww}
\bysame, \emph{The algebra of $sl_3(\mathbb{C})$ conformal blocks},
  arXiv:1206.2535 (2012).

\bibitem[Mic11]{mateusz}
Mateusz Micha{\l}ek, \emph{Geometry of phylogenetic group-based models},
  Journal of Algebra \textbf{339} (2011), 339--356.

\bibitem[PS05]{PachSturm}
Lior Pachter and Bernd Sturmfels (eds.), \emph{Algebraic statistics for
  computational biology}, Cambridge University Press, New York, 2005.

\bibitem[Rai12]{Raicu}
Claudiu Raicu, \emph{Secant varieties of segre--veronese varieties}, Algebra
  and Number Theory \textbf{6} (2012), no.~8, 1817--1868.

\bibitem[SS03]{Ph}
Charles Semple and Mike Steel, \emph{Phylogenetics}, Oxford University Press,
  2003.

\bibitem[SS05]{SS}
Bernd Sturmfels and Seth Sullivant, \emph{Toric ideals of phylogenetic
  invariants}, J. Comput. Biology \textbf{12} (2005), 204--228.

\bibitem[Stu96]{Stks}
Bernd Sturmfels, \emph{Gr\"obner bases and convex polytopes}, University
  Lecture Series, vol.~8, American Mathematical Society, 1996.

\bibitem[Sul07]{Sethtfp}
Seth Sullivant, \emph{Toric fiber products}, Journal of Algebra \textbf{316}
  (2007), no.~2, 560 -- 577, Computational Algebra.

\bibitem[SX10]{SX}
Bernd Sturmfels and Zhiqiang Xu, \emph{Sagbi bases of {Cox}-{Nagata} rings},
  Journal of the European Mathematical Society \textbf{12} (2010), 429--459.

\end{thebibliography}
\vskip 1cm
\footnotesize{
Mateusz Micha\l ek\\
\texttt{wajcha2@poczta.onet.pl}\\
Mathematical Institute of the Polish Academy of Sciences, \'{S}w. Tomasza 30, 31-027 Krak\'{o}w, Poland\\
Max Planck Institute for Mathematics, Vivatsgasse 7, 53111 Bonn, Germany\\

}
\end{document}